\documentclass{article}

\usepackage{arxiv}

\usepackage[utf8]{inputenc} 
\usepackage[T1]{fontenc}    
\usepackage{hyperref}       
\usepackage{url}            
\usepackage{booktabs}       
\usepackage{amsfonts}       
\usepackage{nicefrac}       
\usepackage{microtype}      
\usepackage{lipsum}
\usepackage{caption}
\usepackage{graphicx}
\usepackage{float}
\usepackage{amsmath}
\usepackage{mathtools}
\usepackage{color}
\usepackage[title]{appendix}
\usepackage[shortlabels]{enumitem}

\title{Simple Numerical Algorithm for Generating Hamiltonian Cycles and Edge Labels on Planar Cubic Maps}

\author{
  Emily Kendall \\
  Department of Physics\\
  The University of Auckland\\
  Private Bag 92019\\
  Auckland\\
  New Zealand \\
  \texttt{eken000@auckland.ac.nz} \\
}

\begin{document}
\maketitle


\begin{abstract}

In this work we present an algorithm with which any arbitrary cubic planar map may be constructed through successive edge insertion while simultaneously constructing a set of proper edge labels and Hamiltonian cycles for each configuration. We present a publicly available Python implementation of this algorithm, and discuss both theoretical and numerical support for its validity, with reference to the well-known Four Colour Theorem.

\end{abstract}


\section{Introduction}

The study of planar maps has been an active field of research for many decades, with famous results such as the four colour theorem generating ongoing interest \cite{4CP}. Not only have the properties of planar maps been studied in the context of pure mathematics \cite{Tutte1964, Schaeffer1997, Bousquet2011}, but they have also been shown to have utility in the domain of theoretical physics, particularly in studies of quantum field theory and quantum gravity. \cite{Bessis1980, Bouttier2002, KAZAKOV1986}.

Given the importance of planar maps in both mathematics and physics, the development of numerical tools facilitating their study is of great importance. While there exists a wide variety of problems to be studied relating to planar maps, we focus here on the issue of Hamiltonicity of cubic planar maps \cite{SKUPIEN2002163, RICHMOND1985141}, which itself is closely related to the problem of edge-colouring (edge-labelling). 

Indeed, previous work has been undertaken to develop algorithms to determine Hamilton cycles within planar cubic maps \cite{Price1978}. Here we extend upon this work by introducing a simple algorithm with which one may build an arbitrary planar cubic map through successive edge insertion, whilst enumerating after every insertion the proper edge-labellings associated with the resulting map, as well as its Hamilton cycles. We also present a simple open-source Python implementation of this algorithm.

The structure of this paper is as follows. In Section \ref{sec:descripton} we outline the fundamental concepts of the algorithm with the aid of simple illustrations. In Section \ref{sec:implementation}, we describe how this algorithm is implemented within a simple, standalone Python code. In Section \ref{sec:verification}, we describe tests undertaken on the Python implementation, and provide a simple example of the outputs generated from a given execution. Finally, in Section \ref{sec:discussion}, we discuss the merits of this algorithm, and the scope for future work.

\section{Algorithm description}\label{sec:descripton}

The algorithm which we use for the construction and labelling of cubic planar maps may be decomposed into  three main steps: 1) initialisation, 2) cycle construction, and 3) edge insertion and re-labelling. We describe each step in detail in the following subsections. The corresponding Python implementation is available at \url{https://github.com/erckendall/Edge_Labelling.git}.

\subsection{Step 1: Initialisation}

\begin{figure}[!ht]
\centering
\includegraphics[scale = 0.4, trim={0cm 14cm 18cm 2cm}]{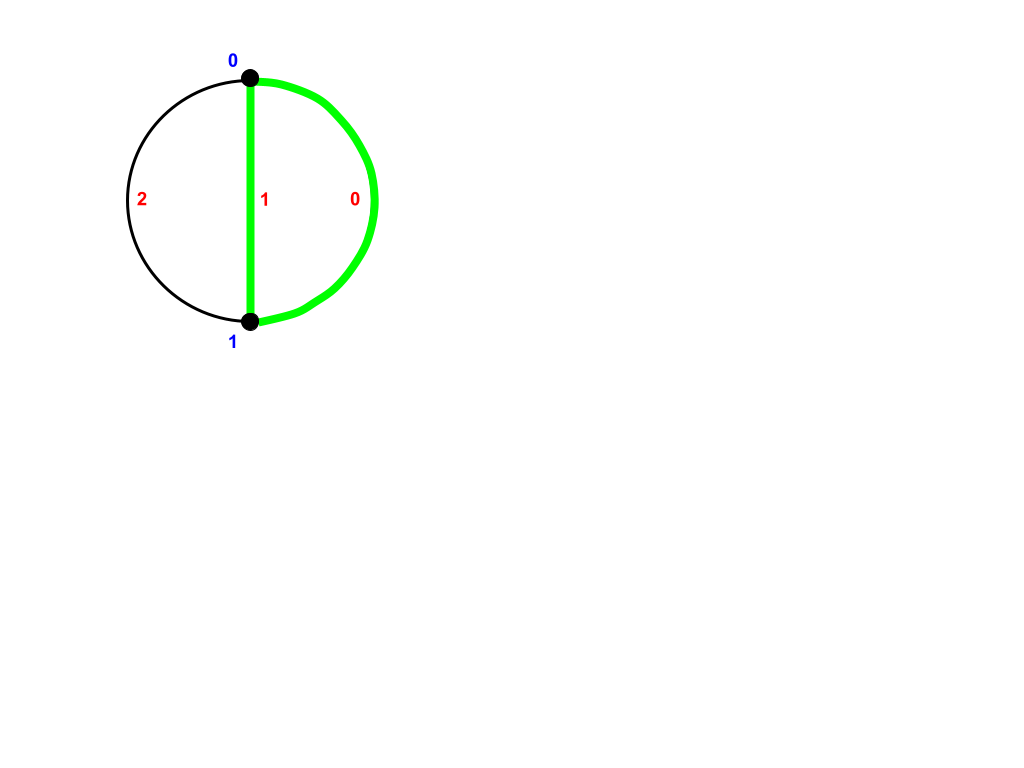}
\caption{The simplest possible initial configuration. A single even cycle through all vertices is shown in green. Edges are enumerated in red, while vertices are enumerated in blue.}\label{fig:simple_config}
\end{figure}

In order to build an arbitrary cubic planar map using this algorithm, we must first define an initial configuration. The initial configuration consists of a planar cubic map and a set of mutually exclusive even cycles over that map, such that every vertex lies on exactly one cycle. The simplest possible initial configuration is the cubic planar map with two vertices, as shown in Figure \ref{fig:simple_config}. Clearly, constructing a set of even cycles through all vertices is trivial in this case. 

\begin{figure}[!ht]
\centering
\includegraphics[scale = 0.4, trim={5cm 9cm 5cm 2cm}]{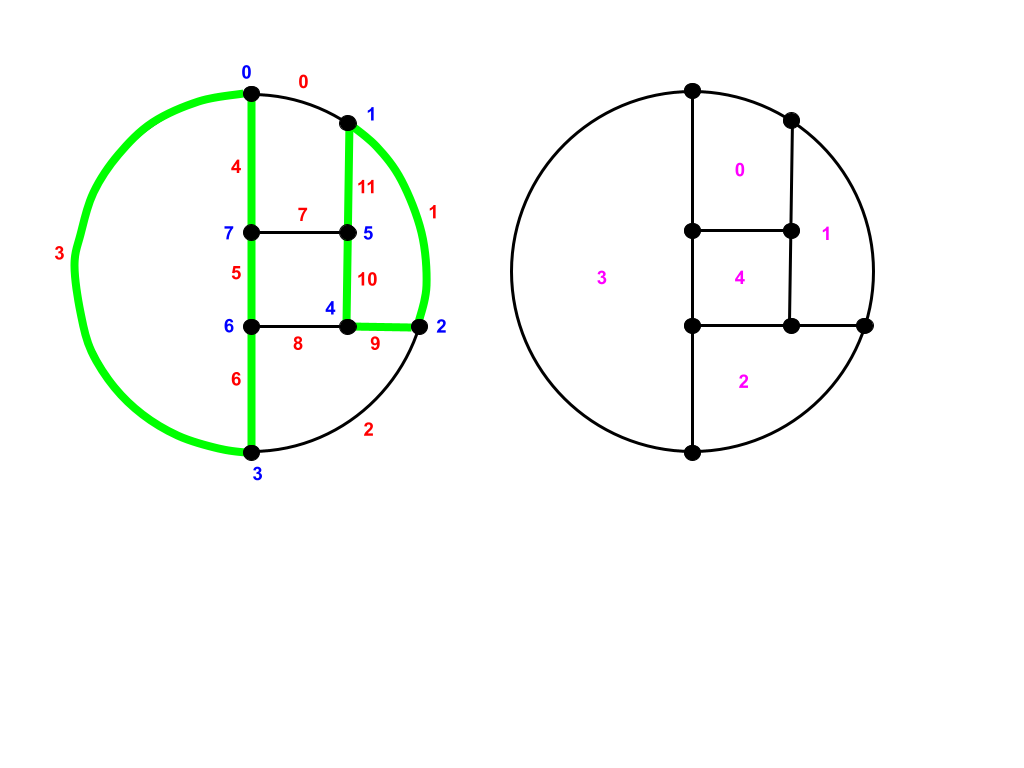}
\caption{An example of a more complicated initial configuration. Here we have a set of two even cycles covering the map (green). Edges are enumerated in red, while vertices are enumerated in blue. We include also an enumeration of the map faces on the right (pink).}\label{fig:initial_config}
\end{figure}

\begin{figure}[!ht]
\centering
\includegraphics[scale = 0.5, trim={0cm 0cm 0cm 0cm}]{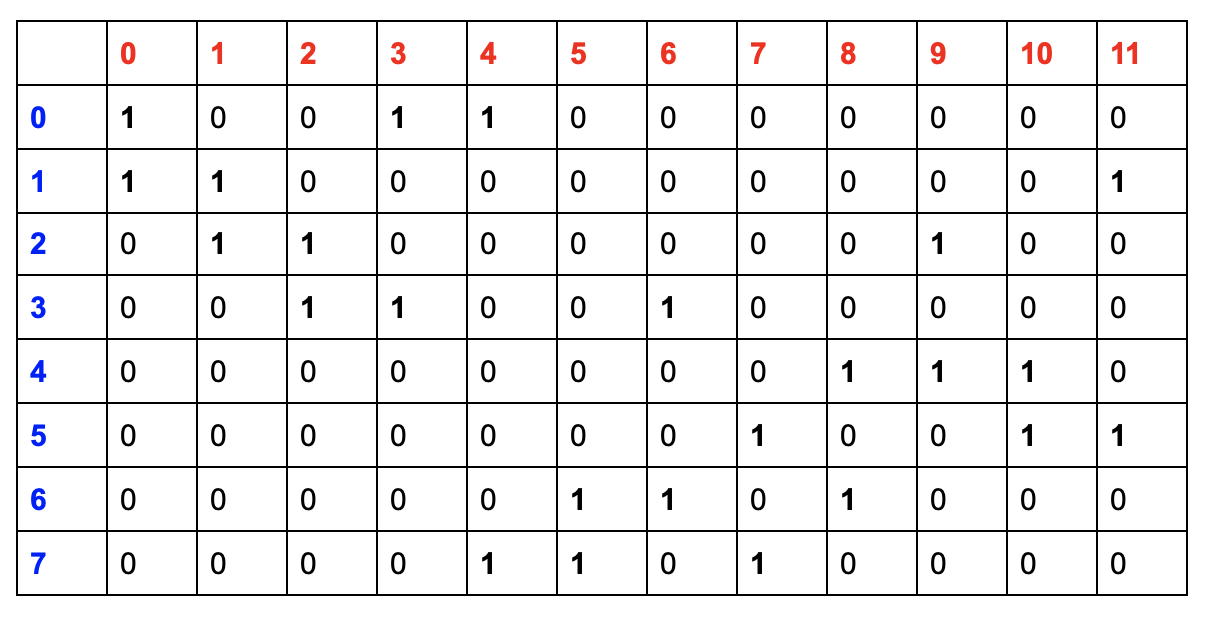}
\caption{Vertex-edge matrix encoding the map illustrated in Figure \ref{fig:initial_config}. Note that for a planar cubic map each column possesses exactly two `1's, while each row possesses exactly three `1's.}\label{fig:vertex_edge}
\end{figure}

\begin{figure}[!ht]
\centering
\includegraphics[scale = 0.5, trim={0cm 0cm 0cm 1cm}]{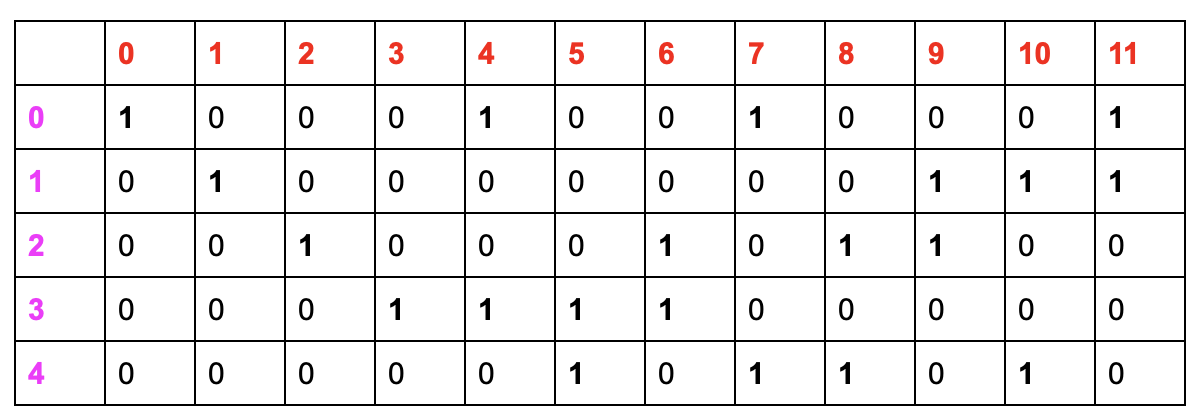}
\caption{Face-edge matrix encoding the map illustrated in Figure \ref{fig:initial_config}. Columns representing external edges possess only one `1', while internal edges possess two. }\label{fig:face_edge}
\end{figure}

While the trivial configuration of Figure \ref{fig:simple_config} may always be used for initialisation, one may also initialise the algorithm using a more complicated map if desired. For the purposes of this explanation, it is instructive to choose the configuration illustrated in Figure \ref{fig:initial_config} as the initial configuration. In this case, we have not one, but two even cycles, each of which traverses four of the eight map vertices. While Figure \ref{fig:initial_config} is useful for illustrative purposes, we must encode this information in matrix form in order to feed it into the Python implementation of the algorithm. To do this, we construct a vertex-edge matrix and a face-edge matrix to describe the map configuration. We also construct a set of lists containing the edges of each cycle. Note that the lists of edges need not be in order, as the Python implementation contains a function to determine this independently. Note also that the numbering of faces, edges and vertices is arbitrary, and requires only that the matrix elements are consistent with the desired graphical representation. For the initial configuration illustrated in Figure \ref{fig:initial_config}, the list of edges which define the cycles is simply $\big[[1,9,10,11], [3,4,5,6]\big]$. Meanwhile, the vertex-edge and face-edge matrices are illustrated in Figures \ref{fig:vertex_edge} and \ref{fig:face_edge}, respectively. 
Once the two matrices and list of edges on each cycle is input into the Python implementation of the algorithm, the initialisation step is complete.

\subsection{Cycle construction}

In the initialisation step of the algorithm, we defined a both a cubic planar map and a set of even cycles upon that map such that every vertex is traversed by exactly one cycle. Henceforth, we will refer to such a set of cycles as a `complete cycle group' (CCG). 

Before we proceed, we note that the initial CCG immediately gives us a proper edge labelling of the map; because each cycle is even, we can apply alternate $a,b,a,b...$ labels around each cycle. Furthermore, because every map vertex is traversed by exactly one cycle, it is therefore also associated with exactly one edge which is not part of any cycle. Hence, the remaining edges may all be labelled $c$, providing a proper edge labelling of the map. We will discuss this labelling procedure further in the following sections.

In the second step of the algorithm, we will use the initial CCG to construct a larger set of alternative CCGs. We describe this procedure below, and then provide a worked example to illustrate each step.
\begin{enumerate}
    \item For each cycle in the initial CCG, we label consecutive edges alternately as $a,b,a,b...$ such that we can define a set of `$a$' edges and a set of `$b$' edges for each cycle.
    \item Choosing either $a$ or $b$ for each cycle, we enumerate all possible sets of choices, amounting to $2^n$ combinations, where $n$ is the number of cycles. For example, if we have two initial cycles in our original CCG, our choices are $\{(a,a), (a,b), (b,a), (b,b)\}$. 
    \item For each of the $2^n$ combinations, we retain (or `turn on') the portions of each cycle which traverse edges with the chosen label, and delete (or `turn off') those portions of the cycle traversing edges labelled with the other letter.
    \item We now create a series of new CCGs by `turning on' each of the edges which were not traversed by any of the cycles in the initial CCG.\footnote{This process necessarily creates even cycles through every vertex of the map, since we are always free to relabel each of the initial cycles independently through $a\leftrightarrow b$, such that the final configuration amounts to a traversal along $\{a-c\}$, or $\{b-c\}$ paths, which are necessarily even when arising from an initial proper edge labelling.} 
    \item We repeat the above steps for each \textit{new} CCG obtained. As we do so, we find that some CCGs which have been obtained previously will emerge again. Eventually, no new configurations will arise, and we will converge onto a finite set of complete cycle groups.
\end{enumerate}

\begin{figure}[!ht]
\centering
\includegraphics[scale = 0.4, trim={5cm 9cm 5cm 2cm}]{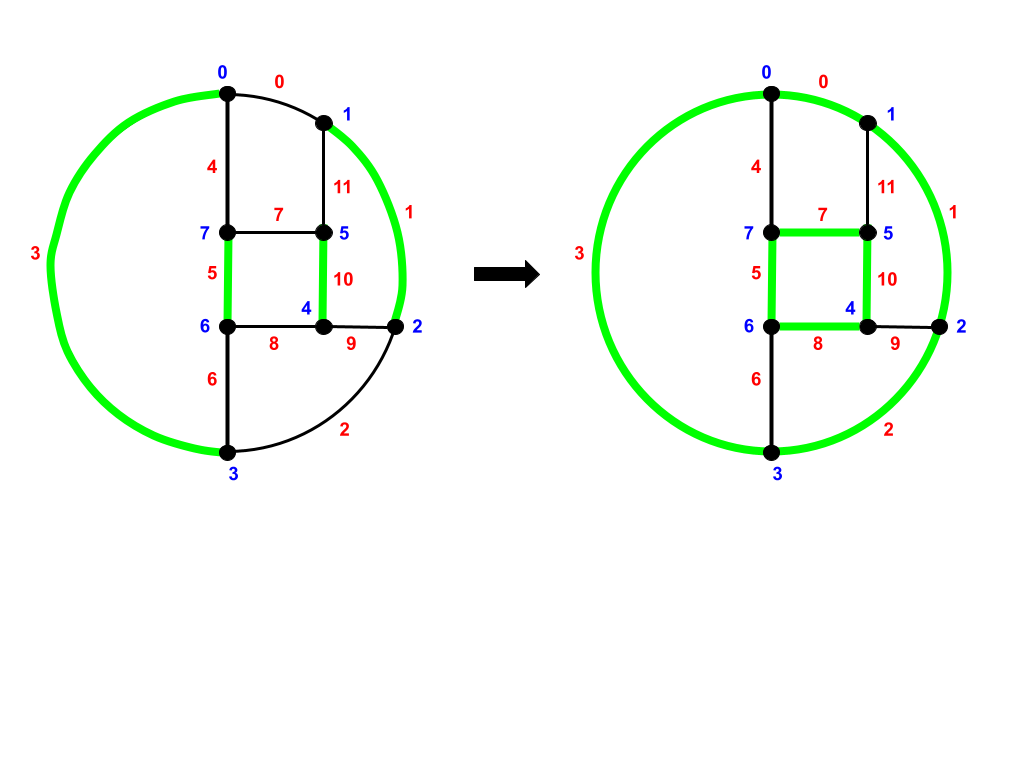}
\caption{Left: first possible combination of original cycle segments from from Figure \ref{fig:initial_config}. Right:
the first new complete cycle group. }\label{fig:set_1}
\end{figure}

\begin{figure}[!ht]
\centering
\includegraphics[scale = 0.4, trim={5cm 9cm 5cm 2cm}]{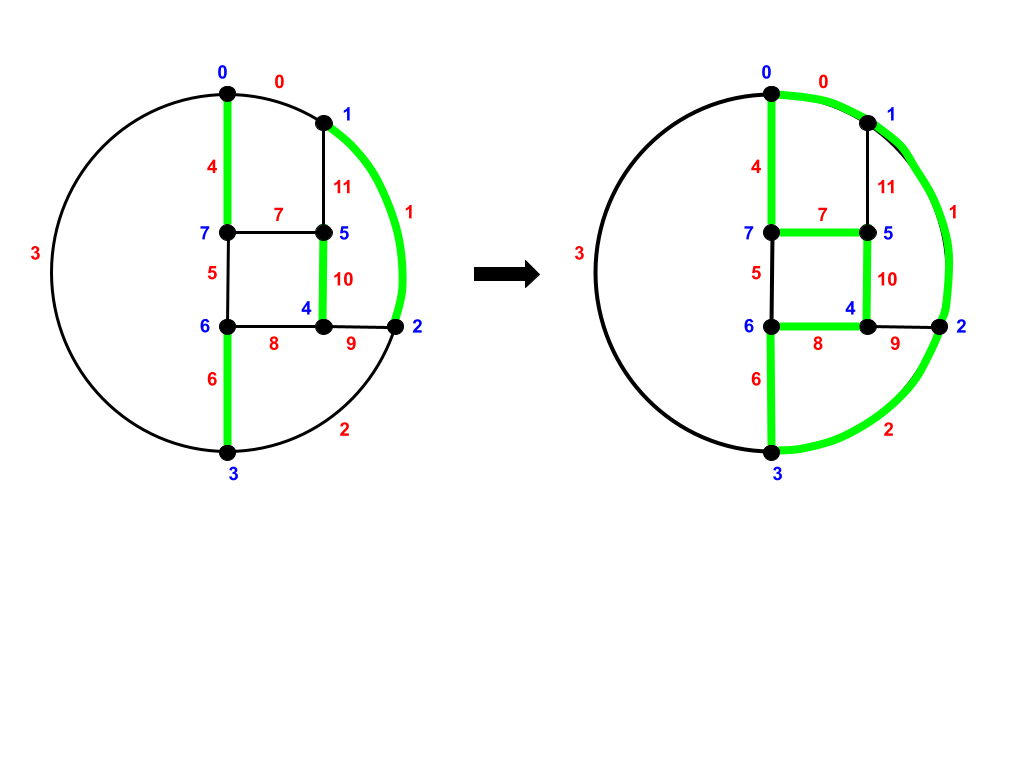}
\caption{Left: second possible combination of original cycle segments from from Figure \ref{fig:initial_config}. Right:
the second new complete cycle group, in this case a single Hamiltonian cycle. }\label{fig:set_2}
\end{figure}

\begin{figure}[!ht]
\centering
\includegraphics[scale = 0.4, trim={5cm 9cm 5cm 2cm}]{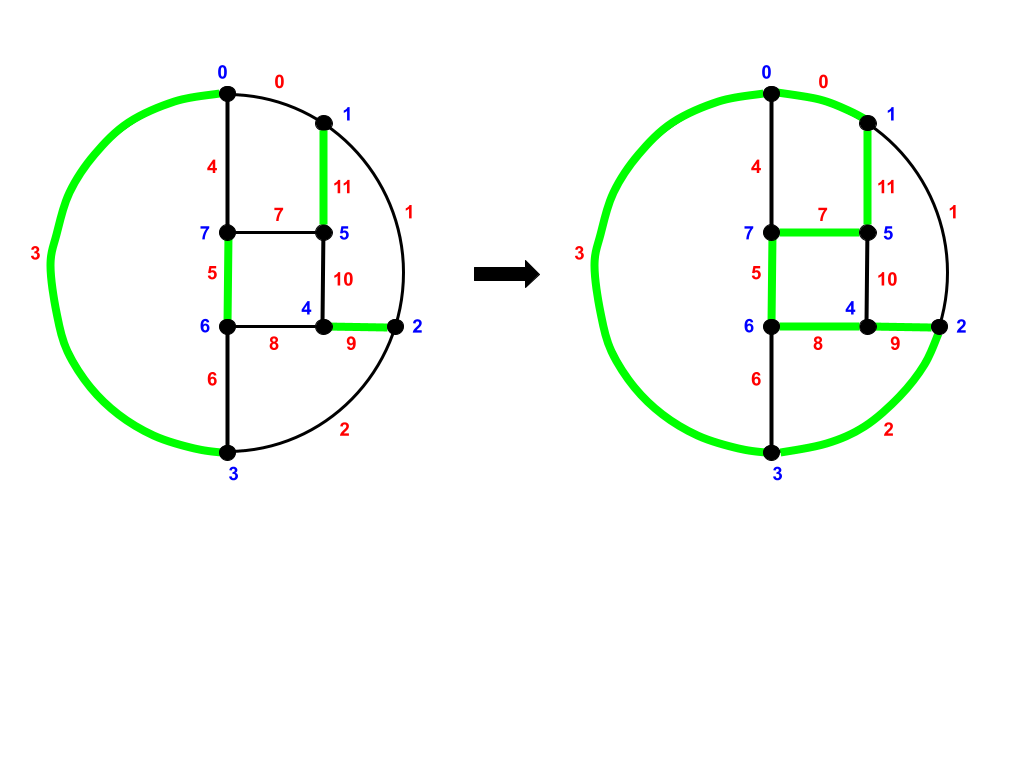}
\caption{Left: third possible combination of original cycle segments from from Figure \ref{fig:initial_config}. Right:
the third new complete cycle group.}\label{fig:set_3}
\end{figure}

\begin{figure}[!ht]
\centering
\includegraphics[scale = 0.4, trim={5cm 9cm 5cm 2cm}]{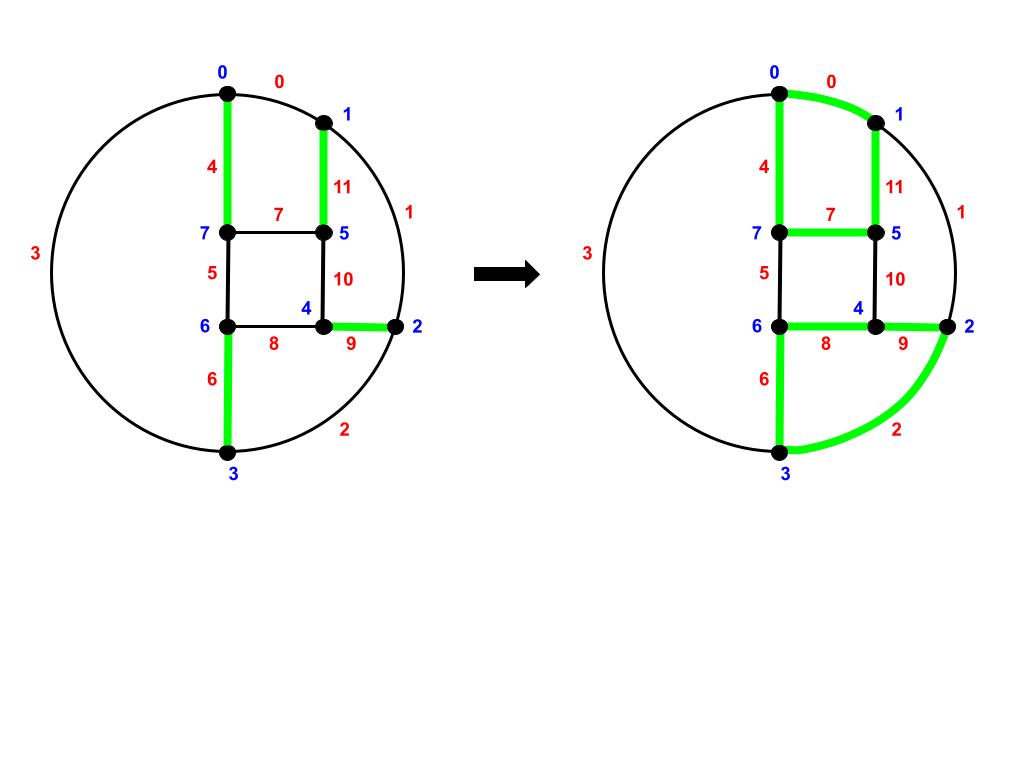}
\caption{Left: fourth and final possible combination of original cycle segments from from Figure \ref{fig:initial_config}. Right:
the fourth new complete cycle group, again a Hamiltonian cycle. }\label{fig:set_4}
\end{figure}

\begin{figure}[!ht]
\centering
\includegraphics[scale = 0.4, trim={0cm 1cm 8cm 2cm}]{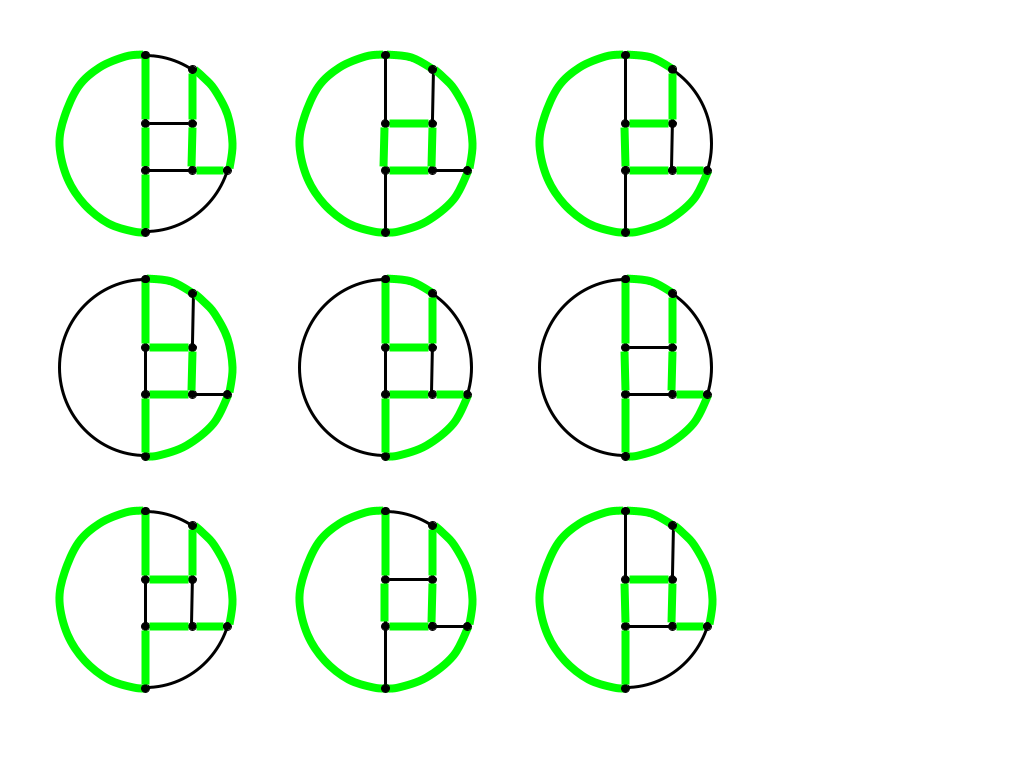}
\caption{The full set of complete cycle groups generated from the initial configuration of Figure \ref{fig:initial_config} using steps 1 to 5 above.}\label{fig:complete_set}
\end{figure}

We now provide an example of this process using the initial configuration of Figure \ref{fig:initial_config}. In this case, we have two cycles in the initial CCG. Therefore, we have four possible combinations of edges to `turn on'. We illustrate these combinations and the resulting CCGs they generate in Figures \ref{fig:set_1} to \ref{fig:set_4}. After having obtained the four new CCGs as illustrated, we then take each of these in turn as our new initial configuration, and repeat the process to produce further CCGs. We do this for every distinct CCG created, until we converge upon a finite set. For this particular example, we converge upon a set of nine distinct CCGs, as illustrated in Figure \ref{fig:complete_set}. We note that of these nine complete cycle groups, six of these are Hamiltonian cycles.

At this point it is pertinent to highlight an important conjecture upon which this algorithm relies:

    \textit{\textbf{Conjecture 1:} If we possess a single CCG for a given planar cubic map, it is possible to produce all possible CCGs for that map through the iterative procedure described in steps 1 to 5, which will always converge to the final complete set of CCGs.}

That is to say, we conjecture that the cycle construction procedure discussed above is exhaustive. We will discuss this notion further in later sections. At present, we simply note that if Conjecture 1 holds, then this algorithm necessarily identifies all possible Hamiltonian cycles for any given map configuration, as all such cycles constitute a CCG.

While the construction of Hamiltonian cycles is in itself interesting, we note that in building the complete set of CCG's, we also build a complete set of proper edge labels, since every proper edge labelling must contain a set of mutually exclusive $a-b-a-b$ cycles, separated by edges which may be labelled $c$. 

While we have already discussed how one may identify a proper edge labelling from any given CCG, it is important to take into account the degeneracy introduced through the arbitrary assignment of `$a$', `$b$', and $`c'$. Hence, we must identify the three distinct sets of edges corresponding to each labelling. Any two labellings are equivalent provided that these sets contain the same elements, irrespective of set order. The Python implementation of this algorithm automatically identifies these equivalences, and outputs only the final set of distinct proper edge labellings.

\subsection{Edge insertion and re-labelling}

In this section we discuss our procedure for introducing new edges into the map while retaining a suitable CCG, and hence a proper edge labelling. This procedure relies on the following important conjecture:

\textit{\textbf{Conjecture 2:} For any two edges which lie on the same face of a planar cubic map, there exists at least one CCG such that both edges lie upon the same cycle within this group.}

\begin{figure}[!ht]
\centering
\includegraphics[scale = 0.4, trim={5cm 9cm 5cm 2cm}]{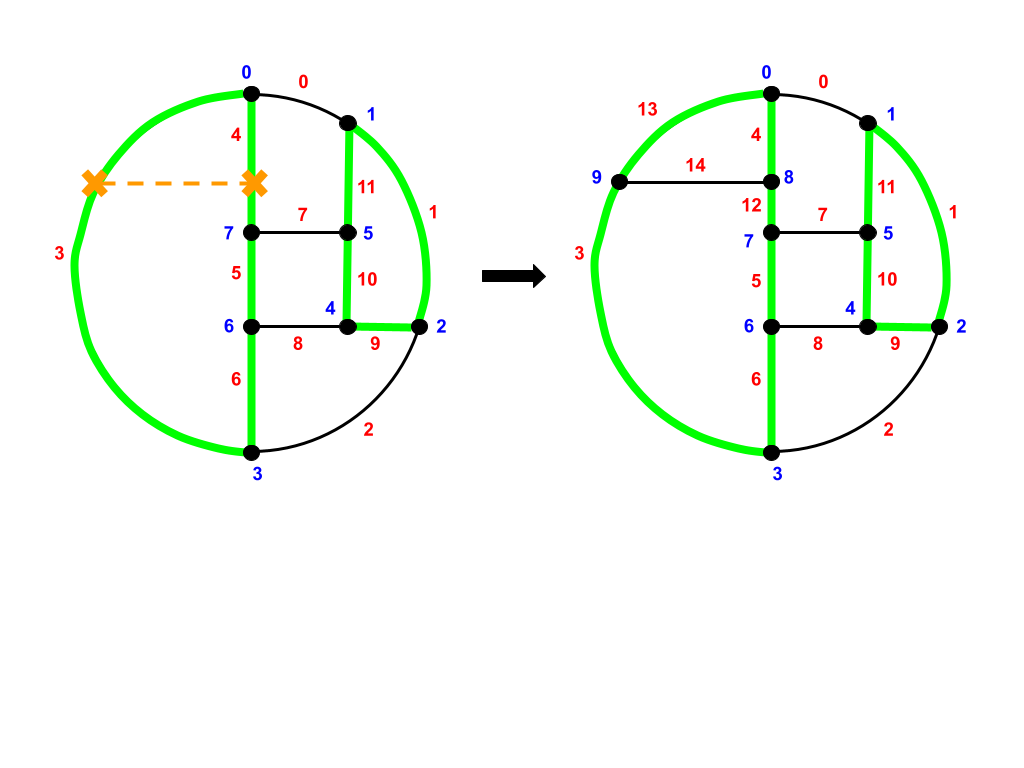}
\caption{Insertion of a new edge between edges 3 and 4 in the original map. From the previously generated set of CCGs, we choose one of these such that both new vertices lie upon the same cycle. We then re-label the map edges and vertices, and obtain a new planar cubic map with a corresponding CCG.}\label{fig:edge_addition}
\end{figure}

We will discuss this conjecture further Section \ref{sec:discussion}. At present, however, we note that if Conjecture 2 holds true, it is always possible to choose a CCG such that the insertion of an edge results in the addition of two new vertices to the \textit{same} cycle within the CCG, thereby preserving evenness, and creating a new planar cubic map for which we automatically have at least one CCG. We illustrate this principle in Figure \ref{fig:edge_addition}.

Having chosen an appropriate CCG, we may then insert the new edge and re-label the map. We then re-define the CCG according to the labelling of the new map, and we thus obtain a new initial configuration to feed into the cycle construction step. Note that we must re-label not only the vertices and edges, but also the faces of the map, as the edge insertion leads to the creation of one additional face. The identification of an appropriate CCG, and the subsequent re-labelling is done automatically in the Python implementation of this algorithm.

\section{Outline of the Python implementation}\label{sec:implementation}

The Python implementation of our algorithm is designed to add edges at random to the initial map configuration until a specified number of edges has been added. After each edge addition, the corresponding CCGs and edge labellings are produced. The Python implementation may be accessed at \url{https://github.com/erckendall/Edge_Labelling.git}.

A skeleton outline of the Python implementation is as follows:

\begin{enumerate}[label*=\arabic*.]
    \item The user inputs the initial vertex-edge matrix (\textit{mat$\_$in}), face-edge matrix (\textit{mat$\_$in$\_$face}), and list of edges on each cycle of the initial CCG (\textit{cycles$\_$in}). The user also specifies the number of edge addition operations to perform using the \textit{iterations} parameter.
    \item Outer loop begins:
    \begin{enumerate}[label*=\arabic*.]
        \item Inner loop begins: 
        \begin{enumerate}[label*=\arabic*.]
            \item The validity of the input vertex-edge matrix is checked using \textit{func$\_$check$\_$validity}.
            \item The edges which are  `turned off' for the current CCG are determined using \textit{func$\_$empties}.
            \item The edges within the list of cycles of the input CCG are ordered using \textit{func$\_$ordering}.
            \item The $2^n$ $a$-$b$ combinations are generated using \textit{func$\_$combos}, where $n$ is the number of cycles in the CCG.
            \item The corresponding set of edge labels for each combination is produced using \textit{func$\_$labels}.
            \item A new list of `turned on' edges is generated for each of the $2^n$ combinations using \textit{func$\_$new$\_$ons}.
            \item The corresponding new CCG for each combination is generated using \textit{func$\_$more$\_$cycles}.
        \end{enumerate}
        \item Inner loop terminates when set of CCGs converges. Output is the complete set of CCG's for the given map and the corresponding edge labellings.
        \item Duplicates in the outputs of the inner loop are removed using \textit{func$\_$duplicate$\_$labels} and \textit{func$\_$remove$\_$duplicates}. Hamiltonian cycles are also identified.
        \item Outputs are printed to the screen for the existing map configuration
        \item One face of the existing map is chosen at random, as are two of its edges (may be the same edge) using \textit{func$\_$choose$\_$edges}.
        \item The order of the edges around the chosen face is determined using \textit{func$\_$order$\_$face}. 
        \item A new vertex-edge matrix is generated corresponding to the addition of a new edge between the two randomly chosen existing edges using \textit{func$\_$new$\_$v$\_$e$\_$mat}
        \item A new face-edge matrix is generated following edge addition using \textit{func$\_$new$\_$f$\_$e$\_$mat} 
        \item An appropriate CCG from the previous map configuration is chosen (such that the two randomly chosen edges both lie on the same cycle) using \textit{func$\_$choose$\_$cycle}.
        \item The chosen CCG is re-written to include the newly created edges using \textit{func$\_$new$\_$cycles}. 
        \item The new matrices and CCG are initialised for the next iteration.
    \end{enumerate}
    \item Outer loop terminates when specified number of edge additions has been completed. Outputs are printed to the screen after each edge addition.
\end{enumerate}

\section{Verification of the Python implementation}\label{sec:verification}

The Python implementation was tested by inserting edges at random and re-labelling the map using the algorithm outlined above. The code can run up to 20 edge additions in under a minute on standard desktop hardware, calculating the set of CCGs and distinct edge labellings after each addition. As this process was randomised and repeated numerous times, we have tested a wide variety of possible map configurations, with the code successfully running to completion in every test.

\begin{figure}[!ht]
\centering
\includegraphics[scale = 0.4, trim={5cm 0cm 5cm 1cm}]{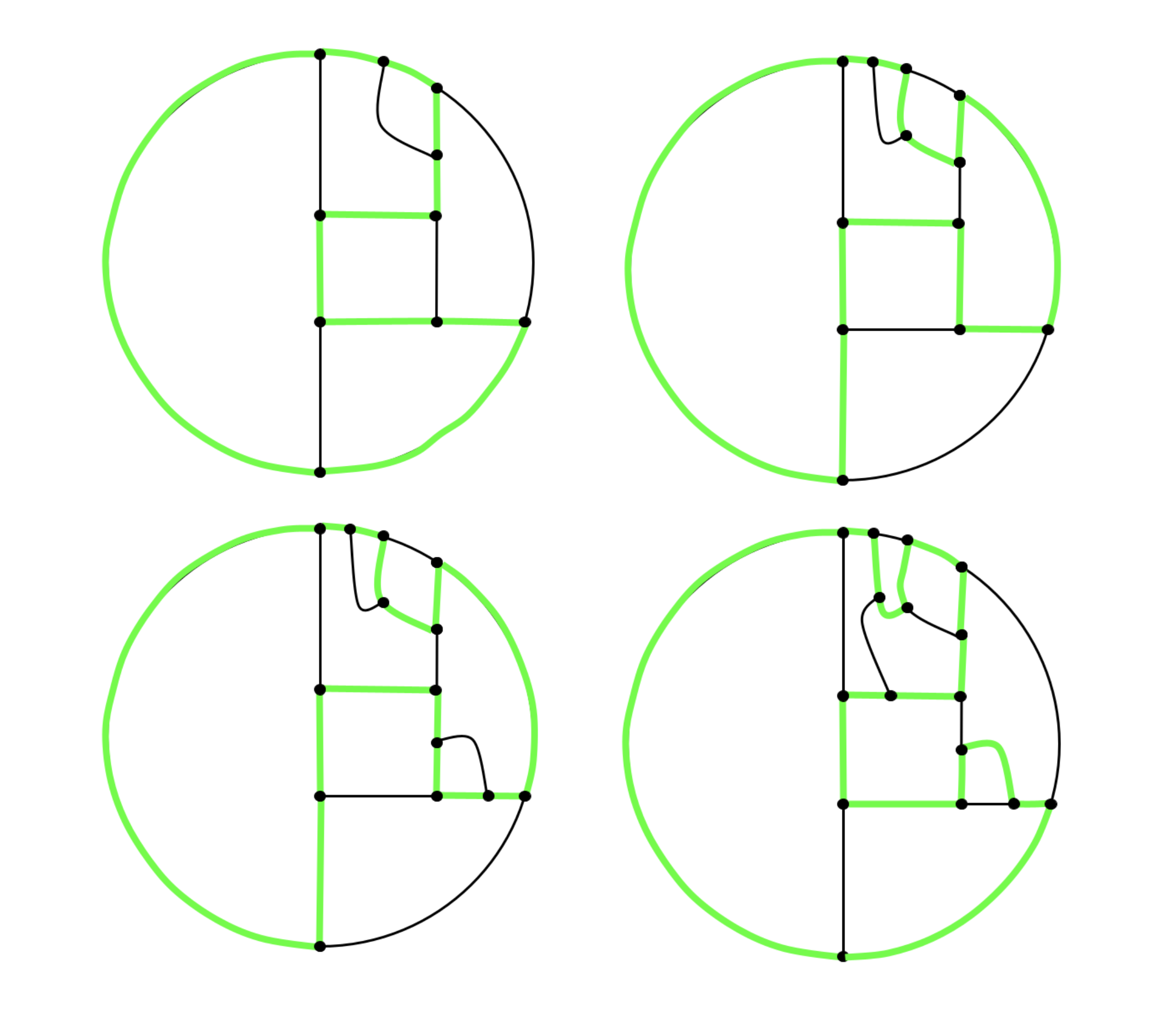}
\caption{Demonstration of random edge addition and CCG generation using the Python implementation. Top left: 1 edge added, top right: 2 edges added, bottom left: 3 edges added, bottom right: 4 edges added.}\label{fig:code_demo}
\end{figure}

Here we demonstrate part of the output of one such run graphically. We take the map illustrated in Figure \ref{fig:initial_config} as our starting configuration, and present here the outcomes of the first four edge additions (Figure \ref{fig:code_demo}), complete with the first CCG assigned in each case. We omit vertex/edge/face labels for clarity.

\section{Discussion and conclusions}\label{sec:discussion}

The problem of edge-labelling of graphs has been extensively studied, (see e.g. \cite{Cranston09}, \cite{BARNOY1992251}, \cite{Biggs72}). The edge-labelling of planar cubic graphs is of particular significance, as it has been established that the four colour theorem (\cite{4CP}) is equivalent to the conjecture that every planar cubic bridgeless graph admits a proper edge labelling \cite{Tait1880} \cite{Chia2008}. We show how a proper edge labelling of a cubic planar map can be mapped into a four colouring in Appendix \ref{sec:labelling}, and we show how this implies a proper four colouring of arbitrary planar maps in Appendix \ref{sec:3_to_4}.

The four colour theorem has been proved through numerical methods (\cite{Appel1977}, \cite{Robertson1996}, \cite{Gonthier2008}), implying therefore that it is always possible to find a proper edge labelling of a cubic planar map. The algorithm discussed here represents a method for determining these proper edge labellings through the construction of even cycle coverings (CCGs) over the map.

The success of this algorithm rests upon the truth of two conjectures, which we reiterate here:

\textit{\textbf{Conjecture 1:} If we possess a single CCG for a given planar cubic map, it is possible to produce all possible CCGs for that map through the iterative procedure described in steps 1 to 5, which will always converge to the final complete set of CCGs.}

\textit{\textbf{Conjecture 2:} For any two edges which lie on the same face of a planar cubic map, there exists at least one CCG such that both edges lie upon the same cycle within this group.}

The first of these conjectures relates to the completeness of the algorithm, while the second may be viewed as a rephrasing of the four colour theorem. That is to say, if Conjecture 2 is true, it is always possible to build an arbitrary planar cubic map possessing a proper edge labelling through successive edge addition, provided the complete set of CCGs can be determined for every intermediate configuration (Conjecture 1). As the four colour theorem has only been proven through numerical means, it is assumed that the proof of Conjecture 2 cannot be reduced to a lesser degree of complexity, but could perhaps also be proved explicitly through computational means. We leave the search for such a proof for future work, but note that the success of the Python implementation of this algorithm presents strong numerical evidence for its validity.

Finally, we note that if Conjecture 1 holds true, then this algorithm presents a method for determining all possible Hamiltonian cycles for a given cubic planar map. The Hamiltonicity of cubic planar maps is also of great interest, see for example \cite{HOLTON1985279}. Our Python implementation automatically checks for the existence of Hamiltonian cycles within the CCGs produced, and will raise an exception if no Hamiltonian cycles are found.\footnote{Note that the Python implementation will also alert the user if none of the Hamiltonian cycles is suitable for the next edge addition, but this is not to say that no Hamiltonian cycles have been found.} We note that throughout our test runs of the Python implementation, we have not encountered a case in which no Hamiltonian cycle is found for any  configuration. 

We intend that this algorithm and its associated implementation may prove a useful tool for the study of the properties of planar cubic maps, and that the conjectures presented here may lead to further progress in the understanding of the implications of the four colour theorem.

\section*{Acknowledgements}

I would like to thank Georgina Liversidge, whose extensive knowledge of graph theory was an invaluable help. I would also like to thank Richard Easther, who patiently allowed me to spend time on graph theory when I was supposed to be doing cosmology.

\bibliographystyle{unsrt}  
\bibliography{references}

\begin{thebibliography}{10}

\bibitem{4CP}
Encyclopedia of~Mathematics.
\newblock Four-colour problem.
\newblock
  \url{http://www.encyclopediaofmath.org/index.php?title=Four-colour_problem&oldid=31833}.
\newblock Accessed: 2019-08-22.

\bibitem{Tutte1964}
W.~G. Brown and W.~T. Tutte.
\newblock On the enumeration of rooted non-separable planar maps.
\newblock {\em Canadian Journal of Mathematics}, 16:572–577, 1964.

\bibitem{Schaeffer1997}
Gilles Schaeffer.
\newblock Bijective census and random generation of eulerian planar maps with
  prescribed vertex degrees.
\newblock {\em The Electronic Journal of Combinatorics [electronic only]},
  4(1):Research paper R20, 14 p.--Research paper R20, 14 p., 1997.

\bibitem{Bousquet2011}
Mireille Bousquet-M{\'e}lou.
\newblock {Counting planar maps, coloured or uncoloured}.
\newblock In {\em {23rd British Combinatorial Conference}}, volume 392 of {\em
  London Math. Soc. Lecture Note Ser}, pages 1--50, Exeter, United Kingdom,
  July 2011.
\newblock Surveys in combinatorics 2011, London Math. Soc. Lecture Note Ser.
  392, pp. 1-50.

\bibitem{Bessis1980}
D.~Bessis, C.~Itzykson, and J.~B. Zuber.
\newblock {Quantum field theory techniques in graphical enumeration}.
\newblock {\em Adv. Appl. Math.}, 1:109--157, 1980.

\bibitem{Bouttier2002}
J.~{Bouttier}, P.~{Di Francesco}, and E.~{Guitter}.
\newblock {Counting colored random triangulations}.
\newblock {\em Nuclear Physics B}, 641(3):519--532, October 2002.

\bibitem{KAZAKOV1986}
V.A. Kazakov.
\newblock Ising model on a dynamical planar random lattice: Exact solution.
\newblock {\em Physics Letters A}, 119(3):140--144, 1986.

\bibitem{SKUPIEN2002163}
Zdzisław Skupień.
\newblock Hamiltonicity of planar cubic multigraphs.
\newblock {\em Discrete Mathematics}, 251(1):163--168, 2002.
\newblock Cycles and Colourings.

\bibitem{RICHMOND1985141}
L.~Bruce Richmond, R.W. Robinson, and N.C. Wormald.
\newblock On hamilton cycles in 3-connected cubic maps.
\newblock In B.R. Alspach and C.D. Godsil, editors, {\em Annals of Discrete
  Mathematics (27): Cycles in Graphs}, volume 115 of {\em North-Holland
  Mathematics Studies}, pages 141--149. North-Holland, 1985.

\bibitem{Price1978}
W.~L. Price.
\newblock An algorithm which generates the hamiltonian circuits of a cubic
  planar map.
\newblock {\em Journal of the London Mathematical Society}, s2-18(2):193--201,
  1978.

\bibitem{Cranston09}
Daniel Cranston.
\newblock Coloring and labeling problems on graphs.
\newblock 02 2009.

\bibitem{BARNOY1992251}
Amotz Bar-Noy, Rajeev Motwani, and Joseph Naor.
\newblock The greedy algorithm is optimal for on-line edge coloring.
\newblock {\em Information Processing Letters}, 44(5):251--253, 1992.

\bibitem{Biggs72}
Norman Biggs.
\newblock An edge-colouring problem.
\newblock {\em The American Mathematical Monthly}, 79(9):1018--1020, 1972.

\bibitem{Tait1880}
P.G. Tait.
\newblock Remarks on the colouring of maps.
\newblock {\em Proc. Roy. Soc. Edinburgh}, 10:501--503, 1880.

\bibitem{Chia2008}
Gek~L. Chia and Siew-Hui Ong.
\newblock Hamilton cycles in cubic graphs.
\newblock 2008.

\bibitem{Appel1977}
Kenneth {Appel} and Wolfgang {Haken}.
\newblock {The Solution of the Four-Color-Map Problem}.
\newblock {\em Scientific American}, 237(4):108--121, Oct 1977.

\bibitem{Robertson1996}
Neil Robertson, Daniel~P. Sanders, Paul Seymour, and Robin Thomas.
\newblock Efficiently four-coloring planar graphs.
\newblock In {\em Proceedings of the Twenty-eighth Annual ACM Symposium on
  Theory of Computing}, STOC '96, pages 571--575, New York, NY, USA, 1996. ACM.

\bibitem{Gonthier2008}
George Gonthier.
\newblock Formal proof - the four color theorem.
\newblock {\em Notices Amer. Math. Soc.}, 55(11), 2008.

\bibitem{HOLTON1985279}
D.A Holton, B~Manvel, and B.D McKay.
\newblock Hamiltonian cycles in cubic 3-connected bipartite planar graphs.
\newblock {\em Journal of Combinatorial Theory, Series B}, 38(3):279--297,
  1985.

\end{thebibliography}

\appendices

\section{Construction of an arbitrary cubic planar map from a trivial starting configuration}\label{sec:construction}

\begin{figure}[!ht]
\centering
\includegraphics[scale = 0.4, trim={2cm 16cm 5cm 4.5cm}]{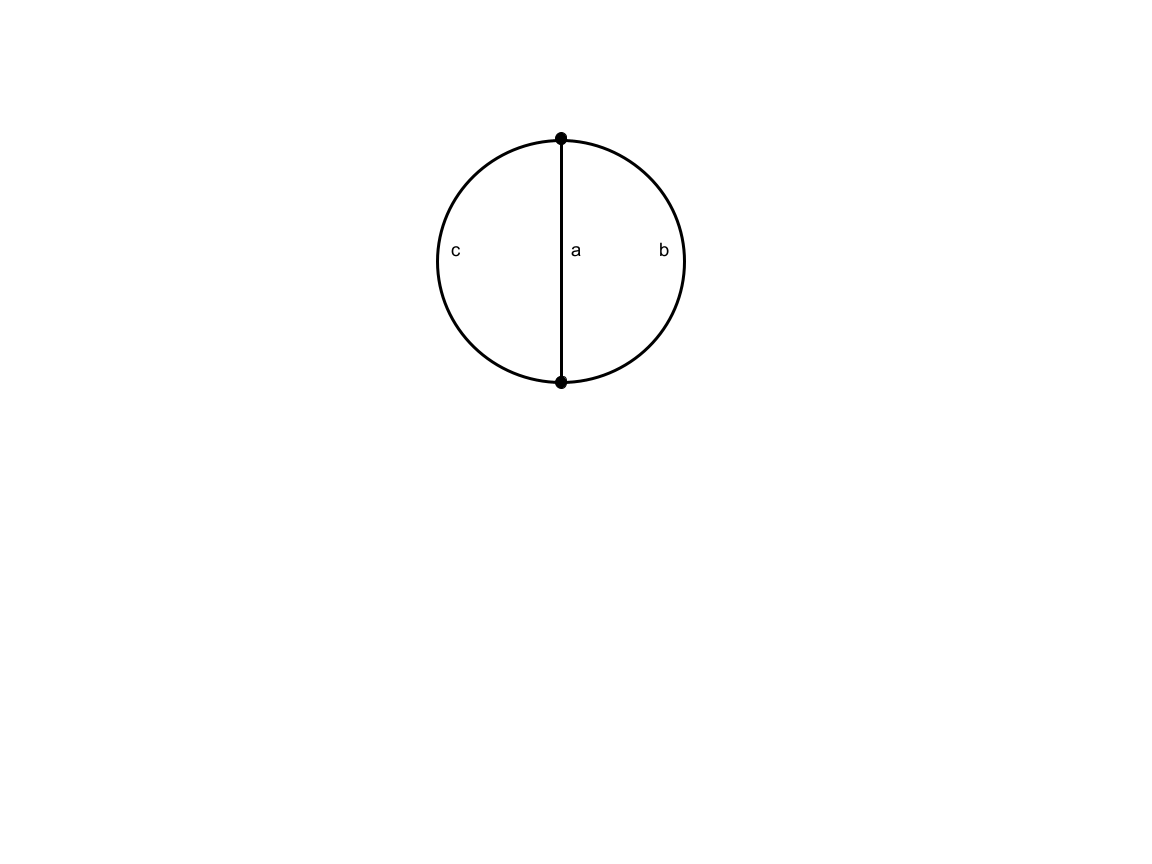}
\caption{The most trivial planar cubic map, from which any arbitrary planar cubic map may be generated through successive insertion of internal edges. Note that this trivial configuration is easily properly edge-labelled, as shown in the figure.}\label{fig:trivial}
\end{figure}

\begin{figure}[!ht]
\centering
\includegraphics[scale = 0.4, trim={2cm 13cm 5cm 3cm}]{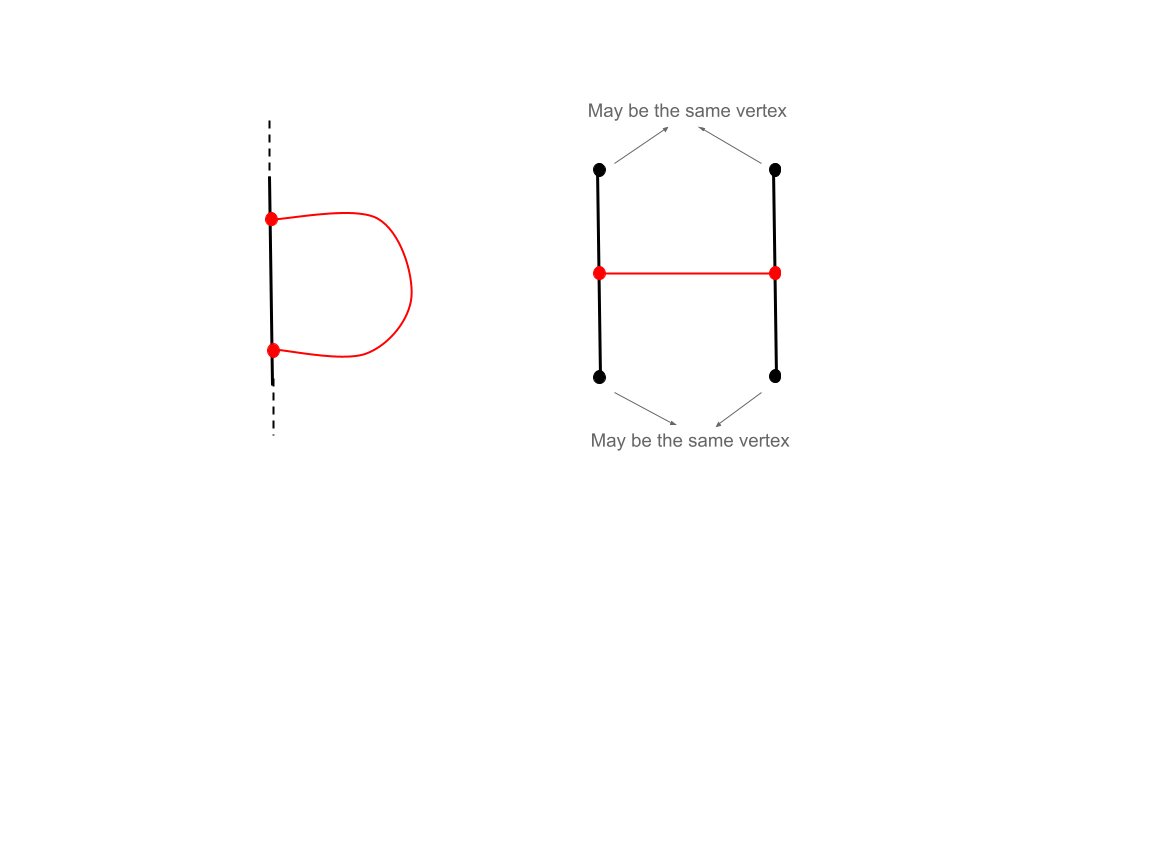}
\caption{Left: We may insert a new edge (red) which terminates along only one existing edge. Right: We may also insert an edge (red) which terminates along two different existing edges (sometimes known as H-insertion).}\label{fig:edge_types}
\end{figure}

All 2-connected planar cubic maps consist of a set of edges which terminate at order-3 vertices only. Hence, if we take an arbitrary cubic planar map and remove any one edge, along with its two terminal vertices, the result is again a 2-connected planar cubic map. 
Because a valid 2-connected map must have a continuous outer boundary, we can consider successively removing only internal edges ( along with their terminal vertices) until we end up with the most trivial planar cubic map, demonstrated in Figure \ref{fig:trivial}. Hence, if every planar cubic map may be reduced to this trivial configuration through successive removal of internal edges, we may conversely construct any arbitrary planar cubic map by successive insertion of edges. A new edge may be inserted such that each of its terminal vertices are placed along the same existing edge, or such that each terminal vertex is placed upon a different existing edge (see Figure \ref{fig:edge_types}). We note, however, that due to planarity the existing edges to which the new edge is attached must always lie upon the same face.

\section{Conversion of arbitrary planar maps into cubic planar maps}\label{sec:3_to_4}

\begin{figure}[!ht]
\centering
\includegraphics[scale = 0.4, trim={10cm 15cm 15cm 3cm}]{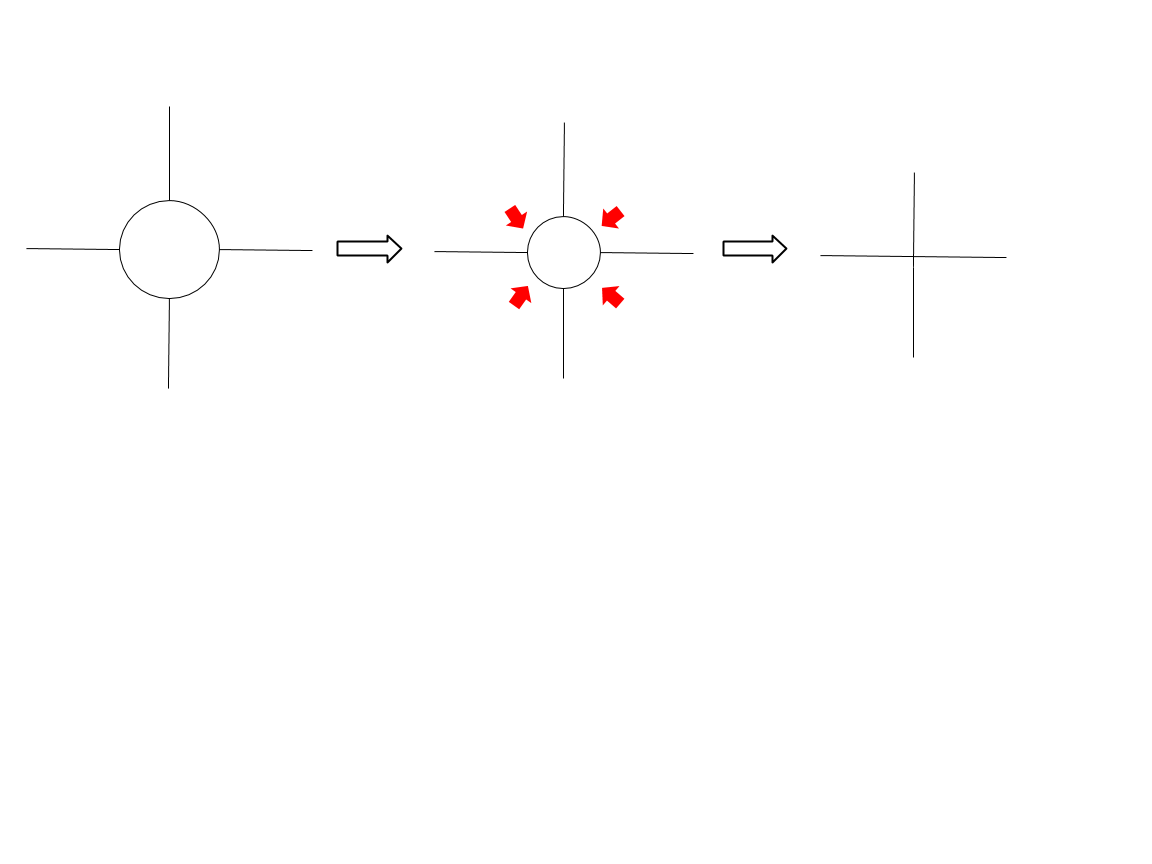}
\caption{Collapsing a circle of order 3 vertices to an infinitesimal point, thus generating a single higher order vertex.}\label{fig:compactify}
\end{figure}

Let us first consider the arrangement of order-3 vertices shown in the left side of Figure \ref{fig:compactify}. here we have a collection of order 3 vertices, arranged around a circle. Let us now imagine constricting the diameter of this circle, as shown in the central diagram. If we continue constricting the circle, we arrive at an infinitesimal point. At this point, the collection of order 3 vertices becomes a single vertex of higher order (in this case order 4). If we consider this process in reverse, we see that we can always transform a vertex of arbitrary order $>3$ into a collection of order 3 vertices. 

Hence, if we take an arbitrary planar map, which contains any number of vertices of order $\geq 3$, we can imagine blowing up each of the vertices in the manner described above, such that we generate a new planar map, containing \textit{only} order 3 vertices. That is, a planar cubic (3-regular) map (Note that the original map must be \textit{bridgeless}, as all edges must constitute borders separating distinct regions. By consequence, no vertices of order <3 exist. Furthermore, any such map must be at least 2-connected, except for `enclaves', which themselves constitute complete 2-connected maps, and can be labelled independently). 

This type of map conversion has the following important consequence: If we can prove that it is always possible to obtain a four-colouring of any 2-connected (bridgeless) planar cubic map, then the full Four Colour Theorem is automatically proved. This is because if we `blow up' an arbitrary planar map in the manner described above, and find that the resulting cubic map possesses a four-colouring, then if we re-collapse the vertices to recover the original map we do not introduce any new neighbours for any region of the original map. Indeed, each region actually loses a neighbour as the circular buffer region collapses into an infinitesimal point, and thus no new restrictions on the colouring of any given region are introduced. Thus, if the expanded (cubic) map is four-colourable, the original map with vertices of arbitrary order must also be four-colourable.

\section{Colouring of a properly edge-labelled 3-connected planar cubic map}\label{sec:labelling}

\begin{figure}[!htb]
\centering
\includegraphics[scale = 0.5, trim={0cm 12cm 17cm 8cm}]{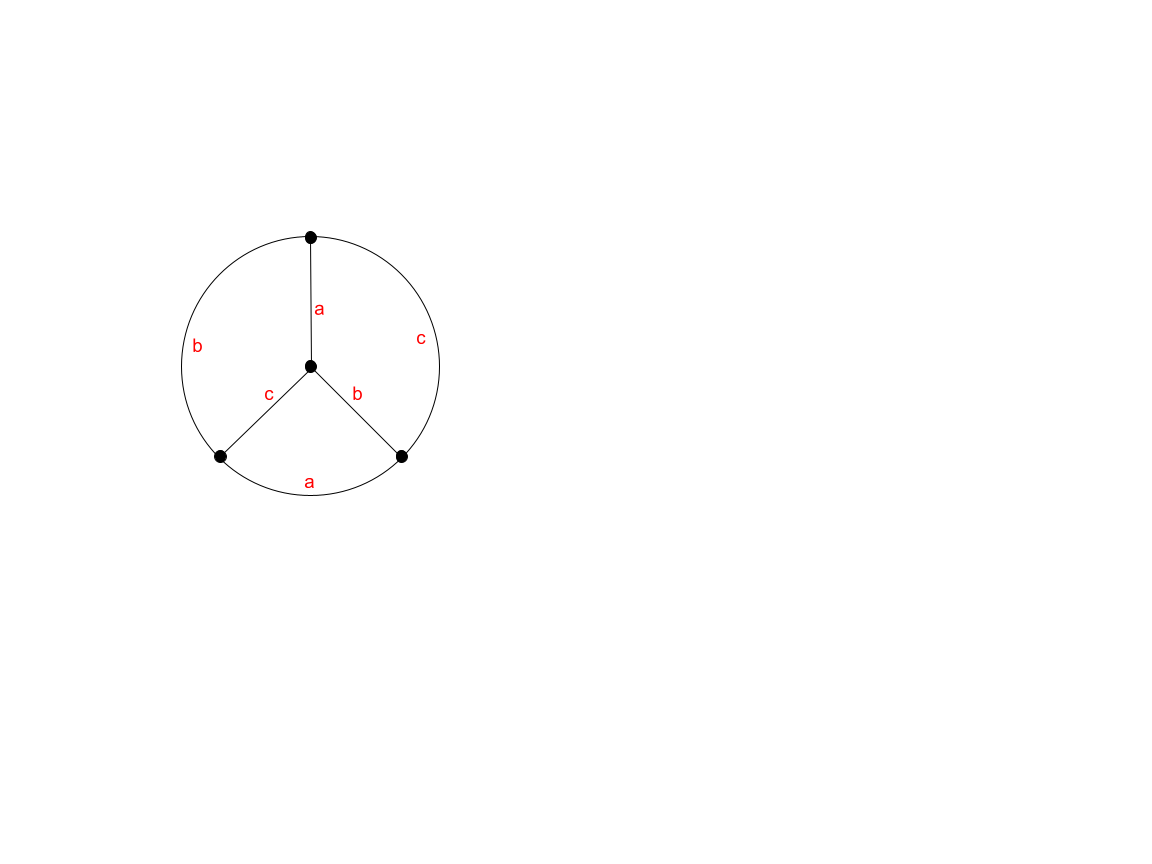}
\caption{Proper edge labelling of the tetrahedral 3-connected planar cubic map.}\label{fig:tetrahedron}
\end{figure}

Let us imagine that we have a properly edge-labelled cubic map, as demonstrated in Figure \ref{fig:tetrahedron}. Here it can be seen that around each vertex, each edge is labelled differently, namely a, b, or c. 


Let us now consider the relationship between the proper edge-labelling and the colouring of the map regions. Assuming we have an ensemble of four different colours to assign to the map regions, for convenience let us denote each colour by a two-index identifier, namely one of the set $\{(+,+), (-,-), (+,-), (-,+)\}$. To satisfy the colouring criterion, no two neighbouring map regions may possess the same index. Hence, when traversing an edge, one of three possible operations must occur:

\begin{enumerate}
    \item The first index swaps sign
    \item The second index swaps sign
    \item Both indices swap sign
\end{enumerate}

\begin{figure}[!ht]
\centering
\includegraphics[scale = 0.5, trim={5cm 19cm 5cm 3cm}]{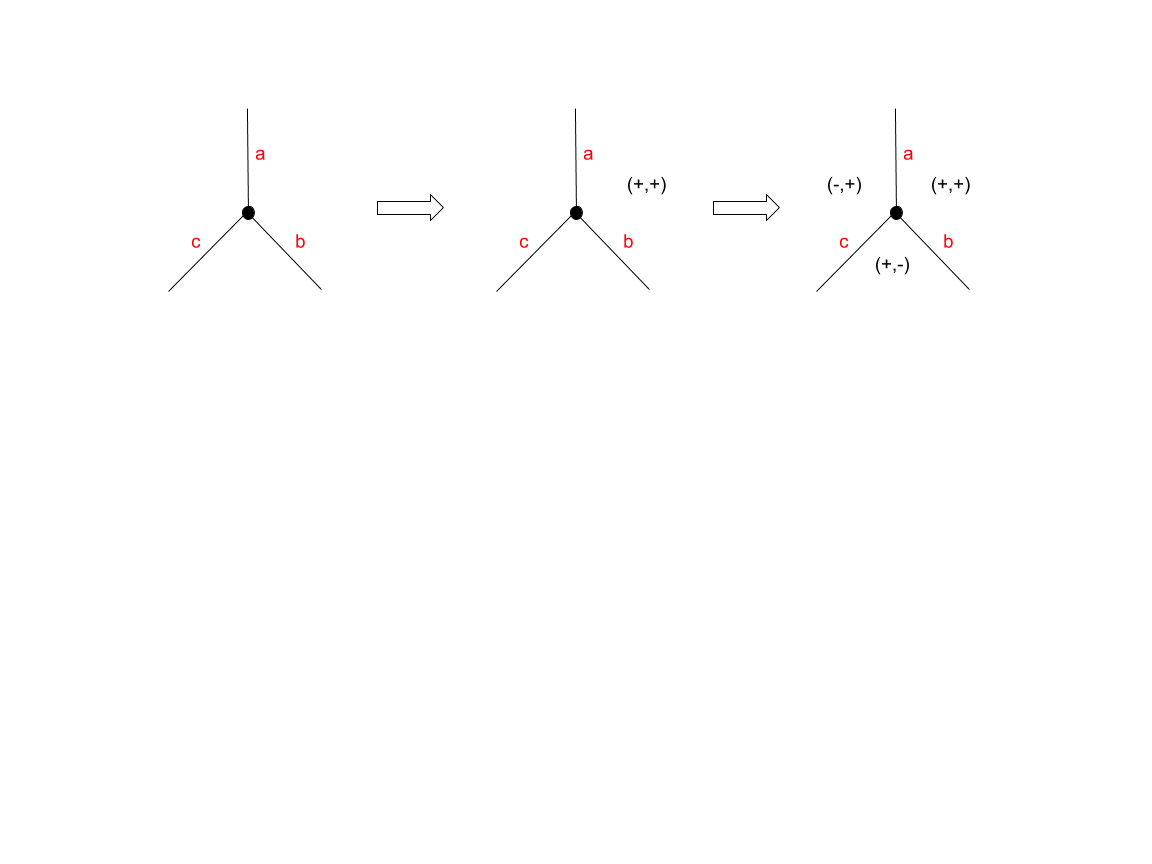}
\caption{Colour index assignment according to edge labelling.}\label{fig:labelling}
\end{figure}

Let us therefore assign to each of the edge labels one of the three above operations. e.g. $a = $ swap first index, $b = $ swap second index, $c = $ swap both indices. Now, if we begin from a properly edge-labelled map such as that shown in Figure \ref{fig:tetrahedron}, we can randomly choose a starting region, and assign to it one of the four two-index identifiers. From here, we can assign identifiers to the rest of the map regions, according to the label of the edge traversed to get to that region. We demonstrate this labelling process in Figure \ref{fig:labelling}.

We can see in the Figure above that because each edge around a vertex is labelled differently, the identifiers of each region around the vertex are mutually consistent. This, in turn, means that as we traverse around the entire border of a given map region vertex by vertex, the identifier assignment will continue to be mutually consistent, such that any path we take to a given region will result in the same identifier. We can be sure that a map labelled in this way will be legally coloured, as an edge must always be traversed between neighbouring regions, changing the identifier. Furthermore, as there are only four possible combinations of the two indices, the map will have no more than four colours. Note also that the "background" upon which the map sits will also acquire a self-consistent colour index, such that the full map + background is in its entirety will be four-colourable. This means that `enclaves', or submaps completely embedded in one larger map region may be treated separately, and inserted without compromising the colouring of the entire map.

\clearpage

\end{document}